\documentclass[12pt,reqno]{amsart}

\usepackage{amssymb, amsthm, amsmath, amsfonts}

\usepackage[marginpar=1cm, margin=2.6cm]{geometry}
\usepackage{verbatim} 
\usepackage{enumerate}

\theoremstyle{remark}
\newtheorem{remark}{Remark}

\theoremstyle{plain}

\title{A note on the divisibility of the Whitehead square}

\author{Haruo Minami}

\address{H. Minami: Professor Emeritus, Nara University of Education}
\email{hminami@camel.plala.or.jp}

\subjclass[2010]{55Q15, 55Q50}

\begin{document}

\begin{abstract}
We show that if we suppose $n\ge 4$ and $\pi_{2n-1}^S$ has no 
2-torsion, then the Whitehead squares of the identity maps on 
$S^{2n+1}$ and $S^{4n+3}$ are divisible by 2. Applying the result of 
G. Wang and Z. Xu on $\pi_{61}^S$, we therefore have that the Kervaire 
invariant one elements in dimensions 62 and 126 exist. 
\end{abstract}

\maketitle

\vspace{-5mm}

\section{Introduction}

Let $[\iota_n, \iota_n]\in \pi_{2n-1}(S^n)$ be the Whitehead square of 
the homotopy class $\iota_n\in \pi_n(S^n)$ of the identity 
map on $S^n$. It is well-known ~\cite{C, G} that this element 
$[\iota_n, \iota_n]$ generates an infinite cyclic subgroup if $n$ is even 
and a cyclic subgroup of order 2 if $n$ is odd $\ne 1, 3, 7$. In particular, 
in the latter case, when $n$ is not of the form $2^r-1$, it splits off as 
a direct summand.  Let $N_k=2^k-1$ $(k \ge 1)$ and put 
$w_k=[\iota_{2N_k+1}, \iota_{2N_k+1}]$. In this note we consider the 
divisibility of $w_k$ by 2. But since we know that $w_1=0$ and $w_2=0$ we assume throughout that $k \ge 3$. 

{\sc Theorem.} {\it Suppose $\pi_{2N_k-1}^S$ has no 2-torsion.   
Then $w_k$ and $w_{k+1}$ are divisible by 2.}

From ~\cite{T} and ~\cite{WX} we know that ${ }_2\pi_{13}^S=0$ and 
${ }_2^{ }\pi_{61}^S=0$ where the subscript $2$ represents the 
2-primary part. Applying these results to the theorem 
we obtain  

{\sc Corollary.} {\it $w_3$, $w_4$, $w_5$ and $w_6$ are divisible by 2.} 

According to ~\cite{BJM}, the existence the Kervaire invariant one 
element $\theta_k\in \pi_{2N_k}^S$ is equivalent to the divisibility of 
$w_k$ by 2. This shows that the corollary assures the existence of 
$\theta_6$ which was unkown. Therefore, by the result of 
 ~\cite{HHR}, we see that the only $\theta_k$ which exist are 
$\theta_1$, $\theta_2$, $\theta_3$, $\theta_4$, $\theta_5$ and 
$\theta_6$. 

In order to prove the theorem, we use an expression for $w_k$ by 
means of the characteristic map of a principal bundle over a sphere. 
Let $T_{n+1}\colon S^{n-1} \to SO(n)$ denote the characteristic 
map of the canonical principal $SO(n)$-bundle $SO(n+1) \to S^n$ 
and let $J$ be the $J$-homomorphism $\pi_{n-1}(SO(n)) \to 
\pi_{2n-1}(S^n)$. Then we know that $[\iota_n, \iota_n]$ 
can be written $[\iota_n, \iota_n]=J([T_{n+1}])$ in particular for $n$ odd ~\cite{GW}, so that $w_k=J([T_{2N_k+2}])$ (the bracket $[ \ ]$ on the right hand side denotes the homotopy class).

Let $\mathbf{R}^n$ be euclidean $n$-space: $x=(x_1, x_2, \ldots, x_n)$.  
Let $S^{n-1}$ be the unit sphere in $\mathbf{R}^n$ with base point 
$x_0=(0, \ldots, 0, 1)$. Then according to  ~\cite{S}, $T_{n+1}$ is  given 
by
\[
T_{n+1}(x)=
\begin{pmatrix}
\delta_{ij}-2x_ix_j 
\end{pmatrix}
\begin{pmatrix}
I_{n-1} &0 \\
0 &-1
\end{pmatrix}
\quad (1 \le i, j \le n)
\]
where $I_{n-1}$ is the identity matrix of dimension $n-1$. 
It is obvious that $T_{n+1}(\pm x_0)=I_n$. Let 
$\Sigma^n=\mathbf{R}^n\cup\{\infty\}$ be the one-point 
compactification of $\mathbf{R}^n$ with $\infty$ as base point. 
Then applying the Hopf construction to $T_{n+1}$ we obtain a map
\[ 
\tau_n\colon \Sigma^n\wedge S^{n-1} \to \Sigma^n         
\]
which satisfies
\begin{equation*}
[\tau_n]=J([T_{n+1}]) \in \pi_{2n-1}(S^n),
\end{equation*}
so that $[\tau_n]=[\iota_n, \iota_n]$ (up to sign) and therefore
\[w_k=[\tau_{2N_k+1}].\]

Let $A_{i_1, \ldots, i_s}$ be the involution on $\mathbf{R}^n$ occurred by 
reversing the signs of the coordinates $x_{i_1}, \ldots, x_{i_s}$ of $x$
where $1\le i_1<\cdots<i_s\le n-1$. Let $\bar{a}_{i_1, \ldots, i_s}$ and $a_{i_1, \ldots, i_s}$ denote the homeomorphisms of $\Sigma^n$ and $S^{n-1}$ onto themselves induced by $A_{i_1, \ldots, i_s}$. Then we have
\begin{equation}
T_{n+1}(a_{i_1, \ldots, i_s}(x))=A_{i_1, \ldots, i_s}T_{n+1}(x)(A_{i_1, \ldots, i_s})^{-1}\qquad (x\in S^{n-1}).
\end{equation}
By applying the definition of $\tau_n$ this yields  
\begin{equation}
\tau_n\circ(\bar{a}_{i_1, \ldots, i_s}\wedge a_{i_1, \ldots, i_s})=
\bar{a}_{i_1, \ldots, i_s}\circ\tau_n.
\end{equation}
Let $\{j_1, \ldots, j_r\}\subset\{i_1, \ldots, i_s\}$ be a subsequence with an even number of elements removed. Then there is a path 
$\mu(t)$ in $SO(n-1)\times\{1\}\subset SO(n)$ joining $A_{j_1, \ldots, j_r}$ and $A_{i_1, \ldots, i_s}$. The path $\mu(t)T_{n+1}(x)\mu(t)^{-1}$ provides a homotopy between $T_{n+1}(a_{j_1, \ldots, j_r}(x))$ and $T_{n+1}(a_{i_1, \ldots, i_s}(x))$ and therefore we see that substituting $\mu(t)$ for $A_{i_1, \ldots, i_{2r}}$ in (1), the equation (2) can be transformed into
\begin{equation}
\tau_n\circ(\bar{a}_{j_1, \ldots, j_r}\wedge a_{i_1, \ldots, i_s})\simeq
\bar{a}_{j_1, \ldots, j_r}\circ\tau_n
\end{equation}
relative to $\Sigma^n\wedge\{x_0, -x_0\}$.

By (2) we have 
$\tau_n\circ(\bar{a}_{1, 2, \ldots, n-1}\wedge a_{1, 2, \ldots, n-1})=
\bar{a}_{1, 2, \ldots, n-1}\circ \tau_n$ and so by (3) 
\begin{equation*}
\begin{split}
\tau_n\circ(1\wedge a_{1, 2, \ldots, n-1})&\simeq\tau_n\qquad\qquad
\qquad\qquad\qquad \ \text{if $n$ is odd},\\
\tau_n\circ(\bar{a}_i\wedge a_{1, 2, \ldots, n-1})&\simeq
\bar{a}_i\circ\tau_n\quad(1\le i\le n-1)\quad \text{if $n$ is even}
\end{split}
\end{equation*}
where $1$ denotes the identity map on $\Sigma^n$. This can be also transformed into the form
\begin{equation}
\begin{split}
\tau_n\circ(1\wedge a_{1, 2})&\simeq\tau_n\qquad\qquad\qquad\qquad
\text{if $n$ is odd},\\
\tau_n\circ(\bar{a}_i\wedge a_i)&\simeq
\bar{a}_i\circ\tau_n\quad(i=1, 3) \quad \ \text{if $n$ is even}
\end{split}
\end{equation}
in a similar way. Below we use this by interpreting it as meaning the above.

Let $D^{n-1}_\pm=\{x\in S^{n-1}\!\mid\! \pm x_1\ge 0\}$ and write
$S^{n-1}_\pm=D^{n-1}_\pm/\partial D^{n-1}_\pm$ for the quotient spaces of them by their boundaries  (in the same sign order). Identify  
$S^{n-1}=D^{n-1}_+\cup D^{n-1}_-$ and  
$S^{n-2}=D^{n-1}_+\cap D^{n-1}_-$ via natural homeomorphisms. Then we have 
$S^{n-1}/S^{n-2}= S^{n-1}_+\vee S^{n-1}_-$ via a natural homeomorphism.
Denote by  
\[\Delta\colon \Sigma^n\wedge S^{n-1} \to (\Sigma^n\wedge S^{n-1}_+)
\vee (\Sigma^n\wedge S^{n-1}_-)\]  the composite of the 
quotient map $\Sigma^n\wedge S^{n-1} \to \Sigma^n\wedge 
S^{n-1}/\Sigma^n\wedge S^{n-2}$ and the canonical homeomorphism 
to $(\Sigma^n\wedge S^{n-1}_+)\vee (\Sigma^n\wedge S^{n-1}_-)$. 
 Let $\pi_\pm\colon S^{n-1} \to S^{n-1}_\pm$ denote the collapsing maps  which are also used to denote $1\wedge \pi_\pm\colon \Sigma^n\wedge S^{n-1}\to \Sigma^n\wedge S_\pm^{n-1}$. 

In general, suppose that given $\mathbf{R}^n$ it is provided with the pair 
$(\Sigma^n, S^{n-1})$ and further that, unless otherwise noted, we denote 
by $\mathbf{R}^{n-i}\subset \mathbf{R}^n$ $(1 \le i \le n)$ the 
subspaces of $\mathbf{R}^n$ spanned by $x$ with the first $i$ 
coordinates equal to zero, and we write $(\Sigma^{n-i}, S^{n-i-1})$ 
for the pairs of their one-point compactifications and unit spheres.
Then by definition of $\tau_n$ we see that 
its restriction to $\Sigma^n\wedge S^{n-i-1}$ can be written as the 
$i$-fold suspension of $\tau_{n-i}$, that is, 
$E^i\tau_{n-i}=\tau_n\!\mid\!\Sigma^n\wedge S^{n-i-1}$.
By definition we also see that the restriction of $\tau_n$ to $\Sigma^{n-i}\wedge S^{n-i-1}$ coincides with $\tau_{n-i}$.  Below we identify 
$\tau_{n-i}=\tau_n\!\mid\!\Sigma^{n-i}\wedge S^{n-i-1}$ for convenience.

It is well known that $\tau_n$  is not null-homotopic except for the case of 
$n=1, 3, 7$. But as a preparation for the proof of Theorem 
we show a result obtained by making an assumption contrary to this 
fact. Below we depend on the use of the method used in the proof 
of the theorem.

{\sc Lemma 1.} {\it Suppose for the sake of argument that 
$\tau_{n-1}\simeq c_\infty$ where $c_\infty$ denotes the constant 
map at $\infty$. Then there exist maps 
\[f_\pm :  \Sigma^n\wedge S^{n-1}_\pm 
\to \Sigma^n
\]
such that 
\[\tau_n\simeq f_+\circ\pi_++f_-\circ\pi_-, \qquad
f_-\circ\pi_- \simeq f_+\circ\pi_+,\]
and therefore $\tau_n\simeq  2f_+\circ\pi_+$.}
\begin{proof} Since $E\tau_{n-1}\simeq c_\infty$ by assumption we see that  there exists a factorization of 
$\tau_n : \Sigma^n\wedge S^{n-1}\to \Sigma^n$ through 
the quotient $\Sigma^n\wedge S^{n-1}/\Sigma^n\wedge S^{n-2}$ into 
the decomposition
\[
\Sigma^n\wedge S^{n-1} \xrightarrow{\Delta}(\Sigma^n\wedge 
S_+^{n-1})\vee (\Sigma^n\wedge S_-^{n-1})\xrightarrow{f_+\vee f_-}
\Sigma^n\vee\Sigma^n\xrightarrow{\mu}\Sigma^n
\]
where $f_\pm : \Sigma^n\wedge S_\pm^{n-1} \to \Sigma^n$ and  
$\mu$ denotes the folding map. Therefore by the definition of the sum of two maps we see that $\tau_n$ decomposes into
\[
\tau_n\simeq f_+\circ \pi_++f_-\circ \pi_-.
\]

By using (4),  observing the behavior of $f_\pm$ on the first two coordinates on $S^n$ we see that 
\begin{equation*}
\begin{split}
(f_-\circ\pi_-)\circ(1\wedge a_{1, 2})&\simeq f_+\circ\pi_+
\qquad\qquad \ \, \text{if $n$ is odd},\\
(f_-\circ\pi_-)\circ(\bar{a}_1\wedge a_1)&
\simeq\bar{a}_1\circ(f_+\circ\pi_+)
\qquad\text{if $n$ is even}.
\end{split}
\end{equation*}
From the first equation it is immediate that, when $n$ is odd,
\[
f_-\circ\pi_-\simeq f_+\circ\pi_+.
\]
For the reason that the assumption is that 
$\tau_{n-1}\simeq c_\infty$, but {\it not its suspension},  it follows that in the case  when $n$ is even, the second equation also represents $f_-\circ\pi_-\simeq f_+\circ\pi_+$ which is the same as in the odd case. Thus we have the lemma.  
\end{proof}

\section{Proof of Theorem}

From now on, let $n=N_k$ (as defined above) and assume 
that the assumption of the theorem is fulfilled, i.e. ${ }_2\pi_{2n-1}^S=0$. 
We also work modulo odd torsion since we consider the 2-primary 
homotopy decomposition of maps. 

From the fact that the suspension homomorphism 
$E\colon \pi_{4n-1}(S^{2n}) \to \pi_{4n}(S^{2n+1})$ of the $EHP$ 
sequence is a surjection with kernel $\mathbf{Z}$ generated by 
$[\iota_{2n}, \iota_{2n}]$ ~\cite{W}, we see that it induces an 
isomorphism ${ }_2\pi_{4n-1}(S^{2n}) \cong{ }_2\pi_{2n-1}^S$ between 
their 2-primary parts. Hence from the assumption above we have 
\begin{equation}\label{e: assump}
\nonumber
\tag{$^\ast$}
_2^{ }\pi_{4n-1}(S^{2n})=0.
\end{equation}

Let $\mathbf{R}^{2n-1}\subset\mathbf{R}^{2n}\subset\mathbf{R}^{2n+1}$; namely let $\mathbf{R}^{2n-1}$ and $\mathbf{R}^{2n}$ be the subspaces of 
$\mathbf{R}^{2n+1}$ consisting of $x=(x_1, \ldots, x_{2n+1})$ with $ x_1=x_2=0$ and $x_1=0$, respectively, and then let
$(\Sigma^{2n-1}, S^{2n-2})\subset (\Sigma^{2n}, S^{2n-1})\subset 
(\Sigma^{2n+1}, S^{2n})$ denote the induced embeddings of their associated 
pairs. Then we have

{\sc Lemma 2.} {\it  In the notation of Lemma 1, under the assumption $(^*)$, there exist maps 
\[
f_\pm\colon  \Sigma^{2n+1}\wedge S^{2n}_\pm \to 
\Sigma^{2n+1}
\]
such that
\[\tau_{2n+1}\simeq f_+\circ\pi_++f_-\circ\pi_-,\qquad  f_-\circ
\pi_-\simeq f_+\circ\pi_+,\]
and therefore $\tau_{2n+1}\simeq 2f_+\circ\pi_+$.}

\begin{proof} 
By  ~\cite{K} we know that $\pi_{2n-2}(SO(2n-1))=\mathbf{Z}_2$ which shows that $2T_{2n}\simeq c_{\mathrm{id}}$ where $c_{\mathrm{id}}$ is the constant map at the identity on $SO(2n-1)$. So by definition of $\tau_{2n-1}$ we have
\[
2\tau_{2n-1}\simeq c_\infty\colon \Sigma^{2n-1}\wedge S^{2n-2} \to \Sigma^{2n-1}.
\] 
Regarding this as the assumption in Lemma 1 in which $\tau_n$ is replaced by 
$2\tau_{2n}$ we can see that, in the notation similar to that used there, it can be written in the form    
\begin{equation*}
\begin{split}
2\tau_{2n}&\simeq f'_+\circ\pi_++f'_-\circ\pi_-\\
&\simeq 2f'_+\circ\pi_+
\end{split}
\end{equation*}
where $f'_\pm\colon\Sigma^{2n}\wedge S_\pm^{2n-1} \to \Sigma^{2n}$.
But due to the assumption $(^*)$ we see that this implies that
\begin{equation}
\tau_{2n}\simeq f'_+\circ\pi_+\colon\Sigma^{2n}\wedge S^{2n-1}\to\Sigma^{2n}.
\end{equation}

This also shows that $\tau_{2n}$ can be homotopically factored through 
\[
\Sigma^{2n}\wedge (S^{2n-1}/S^{2n-2})=(\Sigma^{2n}\wedge S^{2n-1}_+)
\vee(\Sigma^{2n}\wedge S^{2n-1}_-)
\]
and so we are allowed to apply the argument in the proof of Lemma 1 to 
$\tau_{2n+1}$. 
In fact, substituting the suspention $E\tau_{2n}\simeq E(f'_+\circ\pi_+)$ of (5) into 
\begin{equation}
\tau_{2n+1}\circ(1\wedge a_{1, 2})\simeq\tau_{2n+1}
\end{equation}
of (4) we have 
\begin{equation}
E\tau_{2n}=\tau_{2n+1}\!\mid\!\Sigma^{2n+1}\wedge S^{2n-1}\simeq c_\infty
\end{equation}
Therefore, next, taking this as the assumption in Lemma 1 with $\tau_n$ replaced by 
$\tau_{2n+1}$, which is in fact possible through the use of (6) in this case, we can obtain the desired result. Thus we have the lemma. 
\end{proof}

In the coordinate notaion above, let $\mathbf{R}^{2n+2}$ be the subspace of 
$\mathbf{R}^{4n+3}$ consisting of 
\[
x=(x_1, 0, \ldots, 0, x_{2n+3},  \ldots, x_{4n+3})
\]
and then let $(\Sigma^{2n+2}, S^{2n+1})$ denote its associated pair 
in the notation above. But in distinction from the above, let 
$\mathbf{\underline{R}}^{2n+2}$ denote the subspace of $\mathbf{R}^{4n+3}$ consisting of 
\[
\underline{x}=(0, x_2, \ldots, x_{2n+2}, 0, \ldots, 0, x_{4n+3}), 
\]
which is clearly isomorphic to $\mathbf{R}^{2n+2}$, and write 
$(\underline{\Sigma}^{2n+2}, \underline{S}^{2n+1})$ for its associated pair. In particular, $S^{2n+1}\cap\underline{S}^{2n+1}=\{x_0, -x_0\}$. 

Using the coordinates of $x\in S^{2n+1}$ and 
$\underline{x}\in\underline{S}^{2n+1}$ above having the same last coordinates   $x_{4n+3}$ we see that the elements of $S^{4n+2}\subset\mathbf{R}^{4n+3}$ can be expressed in the form 
\[
(rx_1, sx_2, \ldots, sx_{2n+2}, rx_{2n+3}, \ldots, rx_{4n+2}, x_{4n+3}),
\]
written $x(r)$, where $(r, s)\in S^1\subset\mathbf{R}^2$ with 
$r, s\ge 0$. In the setting above we have 

{\sc Lemma 3.} {\it In the notation of Lemma 1, under the assumption $(^*)$, there exist maps 
\[
f_\pm\colon  \Sigma^{4n+3}\wedge S^{4n+2}_\pm \to \Sigma^{4n+3}
\] 
such that
\[\tau_{4n+3}\simeq f_+\circ\pi_++f_-\circ\pi_-, 
\qquad f_-\circ\pi_-\simeq f_+\circ\pi_+,\]
and therefore $\tau_{4n+3}\simeq 2f_+\circ\pi_+$.}

\begin{proof} 
First we prove 
\begin{equation}
\tau_{4n+3}\!\mid\!\Sigma^{2n+2}\wedge S^{2n+1}\simeq c_\infty\quad
\text{and}\quad \tau_{4n+3}\!\mid\! \underline{\Sigma}^{2n+2}\wedge 
\underline{S}^{2n+1}\simeq c_\infty.
\end{equation}
The case of the second equation is quite similar to that of the first one, so we consider only the first case. For simplicity, regard $\mathbf{R}^{2n+2}$ above as the subspace of $\mathbf{R}^{4n+3}$ consisting of 
$x=(x_1, \ldots, x_{2n+2}, 0, \ldots, 0)$ and put $\tau_{2n+2}=\tau_{4n+3}\!\mid\!\Sigma^{2n+2}\wedge S^{2n+1}$. 

Consider its restriction to $\tau_{4n+3}\!\mid\!\Sigma^{2n+2}\wedge S^{2n}$, 
then from the result of Lemma 2, in the notation of maps there, we have 
\[
E\tau_{2n+1}\simeq E(f_+\circ\pi_+)+E(f_-\circ\pi_-).
\]
In the present case, applying this, instead of (5), to
\begin{equation}
\tau_{2n+2}\circ(\bar{a}_3\wedge  a_3)\simeq \bar{a}_3\circ\tau_{2n+2}
\end{equation}
of (4) we can obtain   
\[
-E(f_-\circ\pi_-)\simeq E(f_+\circ\pi_+)
\]
Because, since the associated equation $f_-\circ\pi_-\simeq f_+\circ\pi_+$ 
is obtained based on the use of (5), the same reasoning as in the proof of Lemma 1  can be applied. Thus we have $E\tau_{2n+1}\simeq c_\infty$. 

Putting this into (9) we have that there exist maps 
$f'_\pm\colon\Sigma^{2n+2}\wedge S^{2n+1}_\pm\to\Sigma^{2n+2}$ 
satisfying $\tau_{2n+2}\simeq f'_+\circ\pi_++f'_-\circ\pi_-$ and 
$f'_-\circ\pi_-\simeq f'_+\circ\pi_+$.
But at this point, returning $\mathbf{R}^{2n+2}$ to its original one and 
taking into account the direction of parameter $r$ of the coordinates of 
$S^{4n+2}$, we see that  the latter equation above  must 
be transformed into $f'_-\circ\pi_-\simeq -f'_+\circ\pi_+$ which shows that 
$\tau_{2n+2}\simeq c_\infty$ and therefore we have (8).   

By~\cite{S} we know that $T_{4n+2}$ factors homotopically through 
$SO(4n)\subset SO(4n+1)$ which implies that $E\tau_{4n+1}\simeq c_\infty$,  
since $\pi_{4n}(SO)=\pi_{4n}(SO(4n+2))=0$.  This shows that
\begin{equation}
\tau_{4n+3}\!\mid\!\Sigma^{4n+3}\wedge S^{4n}\simeq c_\infty
\end{equation}
where $S^{4n}=\{x(r)\in S^{4n+2}\!\mid\! x_1=x_2=0\}$. Using this together with the null homotopies of (8) we are able to apply the argument in the proof of Lemma 1 to $\tau_{4n+3}$. In fact, first we see that $\tau_{4n+3}$ factors homotopically through $\Sigma^{4n+3}\wedge Q^{4n+2}$ where 
\[
Q^{4n+2}=S^{4n+2}/(S^{2n+1}\cup\underline{S}^{2n+1}\cup S^{4n}).
\] 

Let $D^{4n+2}_{\mu\nu}=\{x(r)\in S^{4n+2}\!\mid\! \mu x_1\ge 0, \nu x_2\ge 0\}$ 
and put
\[
S^{4n+2}_{\mu\nu}=D^{4n+2}_{\mu\nu}/\partial D^{4n+2}_{\mu\nu} 
\] 
where $\partial D^{4n+2}_{\mu\nu}\subset D_{\mu\nu}^{4n+2}$ is the subspace 
consisting of $x(0), x(1)$ and $x(r)$ with $x_1=x_2=0$. Then we see that 
$S^{4n+2}_{\mu\nu}\approx S^{4n+2}$ and further that $Q^{4n+2}$ can be homeomorphically decomposed into the form
\[
Q^{4n+2}= S^{4n+2}_{++}\vee S^{4n+2}_{--}\vee S^{4n+2}_{+-}
\vee S^{4n+2}_{-+}.
\]
With this identification we consider the homotopy factorization of $\tau_{4n+3}$ above. 

Repeating in mind the process of formation of the null homotopies of (8) and of (10)  above  we see that, in the notation of maps similar to that in the previous  cases, the homotopy factorization of $\tau_{4n+3}$ given above can be written in the form  
\begin{equation}
\tau_{4n+3}\simeq (f_{++}\circ \pi_{++}+f_{--}\circ \pi_{--})
+(f_{+-}\circ \pi_{+-}+f_{-+}\circ \pi_{-+})
\end{equation}
where $f_{\mu\nu}\colon \Sigma^{4n+3}\wedge S_{\lambda\mu}^{4n+2} \to \Sigma^{4n+3}$  and $\pi_{\mu\nu}$ denote the (4n+3)-fold suspensions of  the canonical projections $S^{4n+2}\to S_{\mu\nu}^{4n+2}$ only in the present case. Then in the same way as above, applying the maps $f_{+-}\circ\pi_{+-}$ and $f_{++}\circ\pi_{--}$ in this order to the homotopy expression 
$\tau_{4n+3}\circ(1\wedge a_{1, 2})\simeq \tau_{4n+3}$ of (4) we can obtain 
\begin{equation*}
-f_{-+}\circ\pi_{-+}\simeq f_{+-}\circ\pi_{+-},\qquad f_{--}\circ\pi_{--}\simeq f_{++}\circ\pi_{++}.
\end{equation*}
Substituting these equations into (11) we have  
\[
\tau_{4n+3}\simeq f_{++}\circ \pi_{++}+f_{--}\circ \pi_{--}\simeq 2f_{++}\circ \pi_{++}.
\]
Here, if we interchange $f_{+-}\circ\pi_{+-}$ and 
$f_{++}\circ\pi_{--}$, then similarly we have 
\[
\tau_{4n+3}\simeq f_{+-}\circ \pi_{+-}+f_{-+}\circ \pi_{-+}\simeq 2f_{+-}\circ \pi_{+-}.
\]
But in either case by properly replacing the signs of maps the equation obtained can be rewritten in the form
\[
\tau_{4n+3}\simeq f_+\circ \pi_++f_-\circ \pi_-\simeq 2f_+\circ \pi_+.
\]
This proves the lemma.
\end{proof}

\begin{proof}[Proof of Theorem]
Theorem follows from Lemmas 2 and 3. In fact the results presented there
indicate the divisibility of $w_k$ and $w_{k+1}$ by 2, respectively. 
\end{proof}

\begin{remark} It seems plausible that in the above, if we replace 
$\tau_{2n+1}=\tau_{2N_k+1}$ and $\tau_{4n+3}=\tau_{2N_{k+1}+1}$ by 
$\tau_{4n+3}$ and $\tau_{8n+7}=\tau_{2N_{k+2}+1}$, respectively, then we would be 
able to see that $\tau_{8n+7}$ is divisible by two. But in fact the argument there is 
ineffective because it requires the homotopy expression of $\tau_{4n+2}$ corresponding to (5) which follows from the assumption $(^\ast)$. 
\end{remark}

\end{document}